\documentclass[10pt, reqno]{amsart}
\usepackage{amsmath, amsthm, amscd, amsfonts, amssymb, graphicx, color}
\usepackage[bookmarksnumbered, colorlinks, plainpages]{hyperref}
\hypersetup{colorlinks=true,linkcolor=red, anchorcolor=green, citecolor=cyan, urlcolor=red, filecolor=magenta, pdftoolbar=true}
\usepackage{mathrsfs}
\newtheorem{theorem}{Theorem}[section]
\newtheorem{lemma}[theorem]{Lemma}

\newtheorem{corollary}[theorem]{Corollary}
\theoremstyle{definition}
\newtheorem{definition}[theorem]{Definition}

\theoremstyle{remark}
\newtheorem{remark}[theorem]{Remark}
\numberwithin{equation}{section}

\begin{document}
\setcounter{page}{1}

\title[$L$-functions of noncommutative tori]{$L$-functions of noncommutative tori}

\author[Nikolaev]
{Igor ~V. ~Nikolaev$^1$}

\address{$^{1}$ Department of Mathematics and Computer Science, St.~John's University, 8000 Utopia Parkway,  
New York,  NY 11439, United States.}
\email{\textcolor[rgb]{0.00,0.00,0.84}{igor.v.nikolaev@gmail.com}}

%\dedicatory{In memory of Ola Bratteli}



%\date{Received:  August 14, 2015; Revised: yyyyyy; Accepted: zzzzzz.}

\begin{abstract}
We introduce an analog of the $L$-function for noncommutative tori.
It is proved that such a function coincides with the Hasse-Weil $L$-function
of  an elliptic curve with complex multiplication.  
As a corollary,  one gets a localization formula for
the noncommutative tori with real multiplication. 
\end{abstract}

\maketitle

%**************************************************************************
\section{Introduction}
%***************************************************************************
The Hasse-Weil function $L(\mathcal{E},s)$ is an elegant and powerful  invariant of 
elliptic curves $\mathcal{E}$  [Silverman 1994]  \cite[Chapter 2, Section 10]{S}.  The function $L(\mathcal{E},s)$
encodes important arithmetic information about the elliptic curve $\mathcal{E}$.   
For instance, the famous Birch and Swinnerton-Dyer Conjecture says that the rank of 
$\mathcal{E}$  is equal to the order of zero of  $L(\mathcal{E},s)$ at the point $s=1$
[Tate  1974] \cite[p.198]{Tat1}. Such a rank is always finite by the Mordell Theorem.  
The Hasse-Weil $L$-functions  are also critical
in  the  Langlands Program [Gelbart 1984]  \cite{Gel1}.

Recall that the Sklyanin algebra $S(\alpha,\beta,\gamma)$ is a free $\mathbb{C}$-algebra on  four generators
$\{x_1, \dots, x_4\}$ satisfying six quadratic relations:
$x_1x_2-x_2x_1 = \alpha(x_3x_4+x_4x_3)$,
\quad $x_1x_2+x_2x_1 = x_3x_4-x_4x_3$,
~$x_1x_3-x_3x_1 = \beta(x_4x_2+x_2x_4)$,
~$x_1x_3+x_3x_1 = x_4x_2-x_2x_4$,
~$x_1x_4-x_4x_1 = \gamma(x_2x_3+x_3x_2)$ and 
$x_1x_4+x_4x_1 = x_2x_3-x_3x_2$,
where $\alpha,\beta,\gamma\in \mathbb{C}$ and  
$\alpha+\beta+\gamma+\alpha\beta\gamma=0$ [Stafford \& van ~den ~Bergh 2001]  \cite[Example 8.5]{StaVdb1}.  The algebra $S(\alpha,\beta,\gamma)$ is 
 twisted homogeneous coordinate ring  of an elliptic curve $\mathcal{E}\subset \mathbb{C}P^3$ 
 given in the  Jacobi form $u^2+v^2+w^2+z^2 =
{1-\alpha\over 1+\beta}v^2+{1+\alpha\over 1-\gamma}w^2+z^2 = 0$;
we refer the reader to [Stafford \& van ~den ~Bergh 2001]  \cite{StaVdb1} for the missing definitions and details.
We consider a self-adjoint representation  $\rho:  S(\alpha,\beta,\gamma)\to \mathscr{B}(\mathscr{H})$,  where $\mathscr{B}(\mathscr{H})$ is the ring of 
bounded linear operators on a Hilbert   space $\mathscr{H}$.   The norm-closure of  $\rho(S(\alpha,\beta,\gamma))$ 
is a  noncommutative torus $\mathcal{A}_{\theta}$,    i.e.   a $C^*$-algebra  generated by two unitary operators $u$ and $v$
satisfying the commutation relation $vu=e^{2\pi i\theta}uv$  for a real constant $\theta$.
The map $\mathcal{E}\mapsto \mathcal{A}_{\theta}$ defines  a functor, $F$,  on the category of elliptic curves,   
such that if  $\mathcal{E}$ and $\mathcal{E}'$ are isomorphic,  then 
$\mathcal{A}_{\theta}$ and $\mathcal{A}_{\theta'}$ are Morita equivalent,  i.e.  $\mathcal{A}_{\theta}\otimes \mathscr{K}\cong
\mathcal{A}_{\theta'}\otimes \mathscr{K}$,  where $\mathscr{K}$ is the $C^*$-algebra of compact operators $\theta$ \cite[Section 1.3]{Nik2}.
Moreover, if $\mathcal{E}_{CM}$ is an elliptic curve with complex multiplication [Silverman 1994]  \cite[pp. 95-96]{S},  then 
$F(\mathcal{E}_{CM})=\mathcal{A}_{RM}$, where  $\mathcal{A}_{RM}$ has real multiplication,
i.e.  $\theta$ is a  quadratic irrationality \cite[Theorem 6.1.2]{Nik2}.

The aim of our note is  an $L$-function for noncommutative tori $\mathcal{A}_{RM}=F(\mathcal{E}_{CM})$, such that:
%*******************************************************************************************************
\begin{equation}\label{eq1.1}
L(\mathcal{A}_{RM},s)\equiv L(\mathcal{E}_{CM},s) \quad\hbox{for all}\quad   s\in\mathbb{C}. 
\end{equation}
%******************************************************************************************************
Denote by  $\mathcal{E}_{CM}(\mathbb{F}_p)$  a   localization  of 
elliptic curve  $\mathcal{E}_{CM}$ at the prime ideal $\mathfrak{P}$ over  a prime number $p$
 [Silverman 1994]  \cite[p.171]{S}. 
 Since the cardinals $|\mathcal{E}_{CM}({\Bbb F}_p)|$ generate
the $L(\mathcal{E}_{CM},s)$,  it is clear that  (\ref{eq1.1}) 
implies   a localization formula  for the  algebra $\mathcal{A}_{RM}$  at  $p$,
see corollary  \ref{cor1.2} and remark \ref{rm0}.   A localization theory  of  non-commutative rings is an 
active area of research.

To define an $L$-function of the $\mathcal{A}_{RM}$,  let $\theta$ be a quadratic irrationality 
corresponding to the $\mathcal{A}_{RM}$. 
Denote by $(a_1,\dots,a_k)$ the minimal period of a continued fraction of $\theta$. 
We shall use the following  matrix:
%****************************************************************************************
 \begin{equation}
 A=
 \left(
 \begin{matrix} a_1 & 1\cr 1 & 0\end{matrix}
 \right)
 \left(
 \begin{matrix} a_2 & 1\cr 1 & 0\end{matrix}
 \right)
 \ldots
 \left(
 \begin{matrix} a_k & 1\cr 1 & 0\end{matrix}
 \right). 
 \end{equation}
 %****************************************************************************
For each prime $p$,   one  gets    an  integer matrix:
%************************************************************
\begin{equation}
 L_p=\left(
 \begin{matrix}tr ~A^{\pi(p)} & p\cr -1 & 0\end{matrix}
 \right),
 \end{equation}
 %*****************************************************************
where  $\pi(n)$ an integer-valued function
defined in  Section 2.2.  
Consider an endomorphism  of  the 
$\mathcal{A}_{RM}$  given by  the action of  $L_p$  on the  generators $u$ and $v$.  
The crossed product $C^*$-algebra  $\mathcal{A}_{RM}\rtimes_{L_p}{\Bbb Z}$
  is Morita equivalent  to the  Cuntz-Krieger algebra  
$\mathcal{O}_{L_p}$ [Blackadar 1986]  \cite[Section 10.11.9]{B}. 
For  $z\in {\Bbb C}$ and $\alpha\in\{-1, 0, 1\}$   
we define a  local  zeta function  of  the  $\mathcal{A}_{RM}$ 
%***********************************************************************
$$
\zeta_p(\mathcal{A}_{RM}, z):=
\exp\left(
\sum_{n=1}^{\infty}{|K_0(\mathcal{O}_{\varepsilon_n})|\over n} ~z^n
\right),
~\varepsilon_n=
\begin{cases}
L_p^n, & \mbox{if}  ~p\nmid tr^2(A)-4  \cr
1-\alpha^n,  & \mbox{if}  ~ p ~| ~tr^2(A)-4,
\end{cases}
$$
%*************************************************************************** 
where  $K_0(\mathcal{O}_{\varepsilon_n})$  is  the $K_0$-group of the 
$C^*$-algebra $\mathcal{O}_{\varepsilon_n}$ [Blackadar 1986]  \cite[Chapter III]{B}. 
 %*********************************************************************************
 \begin{definition}
 By an  $L$-function of the  $\mathcal{A}_{RM}$ we understand the  product of the local 
 zetas  taken over all primes:
%********************************************************************************
\begin{equation}
L(\mathcal{A}_{RM}, s)=\prod_{p}  \zeta_p(\mathcal{A}_{RM}, p^{-s}),  \quad s\in {\Bbb C}.
\end{equation}
%********************************************************************************* 
\end{definition}
%**************************************************************************
%***************************************************************************
\begin{theorem}\label{thm1}
$L(\mathcal{A}_{RM}, s)\equiv L(\mathcal{E}_{CM}, s)$.
\end{theorem}
%************************************************************************** 

The following corollary is a localization formula for the algebra $\mathcal{A}_{RM}$.
%***************************************************************************
\begin{corollary}\label{cor1}\label{cor1.2}
$\mathcal{E}_{CM}({\Bbb F}_{p^n})\cong K_0(\mathcal{O}_{\varepsilon_n})$.
\end{corollary}
%************************************************************************** 
\begin{proof}
It is proved in Section 3, that  $\zeta_p(\mathcal{E}_{CM}, z)\equiv\zeta_p(\mathcal{A}_{RM}, z)$
for all $p$.  Hence $|\mathcal{E}_{CM}({\Bbb F}_{p^n})|=|K_0(\mathcal{O}_{\varepsilon_n})|$.
We leave it to the reader to prove that the last equality implies an isomorphism of two finite
abelian groups.
\end{proof}
%**************************************************************************
\begin{remark}\label{rm0}
Corollary \ref{cor1.2} says that  
 the crossed product $\mathcal{A}_{RM}\rtimes_{L_p}{\Bbb Z}$
 is a non-commutative analog of  localization  of  the  polynomial rings 
 at   a   prime ideal $\mathfrak{P}$. 
\end{remark}
%***********************************************************************
The structure of the article is as follows.  
In Section 2,  we introduce notation and review preliminary facts. 
 Theorem \ref{thm1} is proved in Section 3.

%**************************************************************************
\section{Preliminaries}
%***************************************************************************
We briefly review complex multiplication, function $\pi(n)$ and 
Cuntz-Krieger algebras.  For a detailed account, we refer the reader to [Silverman 1994]  \cite{S},
 [Hasse  1950]  \cite{H} and [Cuntz  \&  Krieger 1980]  \cite{CuKr1}, respectively.

%**************************************************************************
\subsection{Complex multiplication}
%***************************************************************************
Let  $\Lambda=\omega_1{\Bbb Z}+\omega_2{\Bbb Z}$ be a lattice in the complex plane ${\Bbb C}$.
Recall that $\Lambda$  defines an elliptic curve $E({\Bbb C}): y^2=4x^3-g_2x-g_3$
via the complex analytic map ${\Bbb C}/\Lambda\to E({\Bbb C})$ given by
the formula $z\mapsto (\wp (z,\Lambda), \wp'(z,\Lambda))$, where  
$g_2=60\sum_{\omega\in\Lambda^{\times}}\omega^{-4}$,
$g_3=140\sum_{\omega\in\Lambda^{\times}}\omega^{-6}$,
$\Lambda^{\times}=\Lambda-\{0\}$ and 
$\wp (z,\Lambda)=z^{-2}+\sum_{\omega\in\Lambda^{\times}} ((z-\omega)^{-2}-\omega^{-2})$
is the Weierstrass $\wp$ function. We shall further identify the elliptic curves $E({\Bbb C})$ 
with the complex tori ${\Bbb C}/\Lambda$. If $\tau=\omega_2/\omega_1$ (a complex modulus),
then $E_{\tau}({\Bbb C}), E_{\tau'}({\Bbb C})$ are isomorphic whenever  
$\tau'\equiv \tau ~mod~SL_2({\Bbb Z})$.  
Recall,  that if $\Lambda$ is a lattice in the complex plane ${\Bbb C}$, then
the endomorphism  ring $End~(\Lambda)$ is isomorphic either to ${\Bbb Z}$ or to
an order, $R$, in the imaginary quadratic number field $k$ [Silverman 1994]  \cite{S}. In the second case,
the lattice is said to have a {\it complex multiplication}. We shall  denote the corresponding
 elliptic  curve by $\mathcal{E}_{CM}$. Consider the cubic $E_{\lambda}: y^2=x(x-1)(x-\lambda)$, $\lambda\in {\Bbb C}-\{0,1\}$.
The $j$-invariant of $E_{\lambda}$ is given by the formula
$j(E_{\lambda})=2^6 (\lambda^2-\lambda+1)^3\lambda^{-2}(\lambda-1)^{-2}$.
To find $\lambda$ corresponding to the $\mathcal{E}_{CM}$,
one has to solve the polynomial equation $j(\mathcal{E}_{CM})=j(E_{\lambda})$ with respect
to $\lambda$. Since $j(\mathcal{E}_{CM})$ is an algebraic integer ([Silverman 1994]  \cite[p.38, Prop.4.5 b]{S}), 
the $\lambda_{CM}\in K$,
where $K$ is an algebraic extension (of degree at most six) of the field
${\Bbb Q}(j(\mathcal{E}_{CM}))$. Thus, each $\mathcal{E}_{CM}$ is isomorphic to a  cubic
$y^2=x(x-1)(x-\lambda_{CM})$ defined over the field $K$.  We shall write this fact as
 $\mathcal{E}_{CM}\cong E(K)$.

Let $K$ be a number field and $E(K)$ an elliptic curve over $K$. 
For each prime ideal $\mathfrak{P}$ of $K$, let ${\Bbb F}_\mathfrak{P}$
be a residue field of $K$ at $\mathfrak{P}$ and 
$q_\mathfrak{P}=N^K_{\Bbb Q}\mathfrak{P}=\# {\Bbb F}_\mathfrak{P}$,
where $N^K_{\Bbb Q}$ is the norm of the ideal $\mathfrak{P}$.
If $E(K)$ has a good reduction at $\mathfrak{P}$, 
one defines
$a_\mathfrak{P}=q_\mathfrak{P}+1-\# \tilde E({\Bbb F}_\mathfrak{P})$,
where $\tilde E$ is a reduction of $E$ modulo the prime ideal
$\mathfrak{P}$. If $E$ has good reduction at $\mathfrak{P}$,
the polynomial
%********************************************************
$
L_\mathfrak{P}(E(K),T)=1-a_\mathfrak{P}T+q_\mathfrak{P}T^2,
$
%*******************************************************
is called the {\it local $L$-series} of $E(K)$ at $\mathfrak{P}$.
If $E$ has bad reduction at $\mathfrak{P}$, the local $L$-series
are $L_\mathfrak{P}(E(K),T)=1-T$ (resp. $L_\mathfrak{P}(E(K),T)=1+T$; $L_\mathfrak{P}(E(K),T)=1$)
if $E$ has split multiplicative reduction at $\mathfrak{P}$ (if $E$ has non-split
multiplicative reduction at $\mathfrak{P}$; if $E$ has additive reduction at $\mathfrak{P}$). 
The global $L$-series defined by the Euler product
%***************************************************
$
L(E(K),s)=\prod_\mathfrak{P} [L_\mathfrak{P}(E(K), q_\mathfrak{P}^{-s})]^{-1},
$
%****************************************************   
is called a {\it Hasse-Weil $L$-function} of the elliptic curve $E(K)$.

Let $A_K^*$ be the idele group of the number field $K$. 
A continuous homomorphism $\psi: A_K^*\to {\Bbb C}^*$ with
the property $\psi(K^*)=1$ is called a {\it Gr\"ossencharacter}
on $K$. (The asterisk denotes the group of invertible elements
of the corresponding ring.) The {\it Hecke $L$-series} attached
to the Gr\"ossencharacter $\psi: A_K^*\to {\Bbb C}^*$ is defined 
by the Euler product
%************************************************************
$
L(s,\psi)=\prod_\mathfrak{P} (1-\psi(\mathfrak{P})q_\mathfrak{P}^{-s})^{-1},
$
%*************************************************************
where the product is taken over all prime ideals of $K$.

Let $\mathcal{E}_{CM}\cong E(K)$ be an elliptic curve with complex multiplication
by the ring of integers $R$ of an imaginary quadratic field $k$, and assume
that $K\supset k$. Let $\mathfrak{P}$ be a prime ideal of $K$ at which $E(K)$
has a good reduction. If  $\tilde E$ is a reduction of $E(K)$ at $\mathfrak{P}$,
we let $\phi_\mathfrak{P}: \tilde E\to\tilde E$ be the associated Fr\"obenius map. 
Finally, let $\psi_{E(K)}: A_K^*\to k^*$ be the Gr\"ossencharacter attached
to the $\mathcal{E}_{CM}$,  see [Silverman 1994]  \cite[p.168]{S}. The following diagram is known 
to be commutative:

%*******************************************************************
\begin{picture}(300,110)(-80,-5)
\put(20,70){\vector(0,-1){35}}
\put(130,70){\vector(0,-1){35}}
\put(45,23){\vector(1,0){53}}
\put(45,83){\vector(1,0){53}}
\put(15,20){$\tilde E$}
%\put(5,55){$F$}
%\put(140,55){$F$}
\put(123,20){$\tilde E$}
\put(15,80){$E(K)$}
\put(115,80){$E(K)$}
\put(60,30){$\phi_\mathfrak{P}$}
\put(50,90){$\psi_{E(K)}(\mathfrak{P})$}
\end{picture}
%***********************************************************

\medskip\noindent
see [Silverman 1994]  \cite[p.174]{S}. In particular, $\psi_{E(K)}(\mathfrak{P})$
is an endomorphism of the $E(K)$ given by the complex number
$\alpha_{E(K)}(\mathfrak{P})\in R$. By $\overline{\psi}_{E(K)}(\mathfrak{P})$
one understand the conjugate Gr\"ossencharacter viewed as a complex number.
The {\it Deuring Theorem} says that the Hasse-Weil $L$-function of the $E(K)$
is related to the Hecke $L$-series of the $\psi_{E(K)}$ by the formula
%****************************************************************** 
$
L(E(K),s)\equiv L(s,\psi_{E(K)})L(s, \overline{\psi}_{E(K)}).
$
%******************************************************************

%**************************************************************************
\subsection{Function $\pi(n)$}
%***************************************************************************
Let $\mathfrak{k}={\Bbb Q}(\sqrt{D})$ be a real quadratic number field
and $O_\mathfrak{k}$ its ring of integers. For rational integer $n\ge 1$
we shall write $\mathfrak{R}_n\subseteq O_\mathfrak{k}$  to denote an
order (i.e. a subring containing  $1$) of $O_\mathfrak{k}$.  The order
$\mathfrak{R}_n$ has a basis $\{1, n\omega\}$, where 
%***********************************************
\begin{equation}
\omega=\begin{cases}{\sqrt{D}+1\over 2} & if ~D\equiv 1 ~ mod~4,\cr
               \sqrt{D} & if ~D\equiv 2,3 ~ mod~4.
\end{cases}               
\end{equation}
%***************************************************
In other words, $\mathfrak{R}_n={\Bbb Z}+(n\omega){\Bbb Z}$. 
It is clear, that $\mathfrak{R}_1=O_\mathfrak{k}$ and the fundamental 
unit of $O_\mathfrak{k}$ we shall denote by $\varepsilon$. 
Each $\mathfrak{R}_n$ has its own fundamental unit, which
we shall write as $\varepsilon_n$;  notice that $\varepsilon_n\ne
\varepsilon$  unless $n= 1$. 

There exists the well-known  formula,  which relates $\varepsilon_n$
to the fundamental unit $\varepsilon$,  see e.g. [Hasse  1950]  \cite[p.297]{H}.
Denote by $\mathfrak{G}_n:=U(O_\mathfrak{k}/nO_\mathfrak{k})$ 
the multiplicative group of invertible elements (units) of 
the residue ring  $O_\mathfrak{k}/nO_\mathfrak{k}$;   clearly, all units 
of $O_\mathfrak{k}$ map (under the natural $mod~n$ homomorphism) to $\mathfrak{G}_n$. 
Likewise  let $\mathfrak{G}_n:=U(\mathfrak{R}_n/n\mathfrak{R}_n)$
be the group of units of the residue ring $\mathfrak{R}_n/n\mathfrak{R}_n$;
it is not hard to prove ([Hasse  1950]  \cite[p.296]{H}), that $\mathfrak{G}_n\cong U({\Bbb Z}/n{\Bbb Z})$
the ``rational'' unit group of the residue ring ${\Bbb Z}/n{\Bbb Z}$. 
Similarly, all units of the order $\mathfrak{R}_n$ map to $\mathfrak{G}_n$.  
Since units of $\mathfrak{R}_n$ are also units of $O_\mathfrak{k}$
(but not vice versa),  $\mathfrak{R}_n$ is a subgroup of $\mathfrak{G}_n$;
in particular,   $|\mathfrak{G}_n|/|\mathfrak{R}_n|$ is an integer number
and $|\mathfrak{G}_n|=\varphi(n)$, where $\varphi(n)$ is the Euler totient function. 
In general,  the following formula is true
%***************************************************************************
\begin{equation}\label{eq16}
{|\mathfrak{G}_n|\over |\mathfrak{R}_n|}=n\prod_{p_i | n}
\left(1-\left({D\over p_i}\right){1\over p_i}\right),
\end{equation}
%**************************************************************************** 
where $\left({D\over p_i}\right)$ is the Legendre symbol, 
see [Hasse  1950]  \cite[p.351]{H}.   We shall write $\pi(n)$ to denote   the least integer number  dividing 
$|\mathfrak{G}_n|/|\mathfrak{R}_n|$ and   such that
$\varepsilon^{\pi(n)}$ is  a unit of $\mathfrak{R}_n$ (i.e. belongs to $\mathfrak{G}_n$);
the following Satz XIII$^{\prime}$ of Hasse's book says that the unit is the fundamental
unit of  order $\mathfrak{R}_n$.  
%***************************************************************
\begin{lemma}\label{lm10}
{\bf ([Hasse  1950]  \cite[p.298]{H})}
$\varepsilon_n=\varepsilon^{\pi(n)}.$
\end{lemma}
%**************************************************************

\bigskip\noindent
In what follows,  we deal with the special case $n=p$  is a prime number;
in this case formula (\ref{eq16}) becomes
%***************************************************************************
\begin{equation}\label{eq17}
{|\mathfrak{G}_p|\over |\mathfrak{R}_p|}=p-\left({D\over p}\right).
\end{equation}
%****************************************************************************   
Notice,  that lemma \ref{lm10} asserts  existence  of the number $\pi(n)$ 
(as one of the divisors of  $|\mathfrak{G}_n|/|\mathfrak{R}_n|$) yet no analytic 
formula for $\pi(n)$ is given;  it would be rather interesting to have such
a formula.

%**************************************************************************
\subsection{Cuntz-Krieger algebras}
%***************************************************************************
A  Cuntz-Krieger algebra, $\mathcal{O}_B$, is the $C^*$-algebra
generated by partial isometries $s_1,\dots, s_n$ that act on a Hilbert space in such 
a way that their support projections $Q_i=s_i^*s_i$ and their  range projections
$P_i=s_is_i^*$ are orthogonal and satisfy the relations
$Q_i=\sum_{j=i}^nb_{ij}P_j$,
for an $n\times n$  matrix $B=(b_{ij})$ consisting of $0$'s 
and $1$'s [Cuntz  \&  Krieger 1980]  \cite{CuKr1}.  The notion is extendable to the 
matrices $B$ with the non-negative integer entries 
[Cuntz  \&  Krieger 1980]  \cite[Remark 2.16]{CuKr1}.
It is known,  that the $C^*$-algebra $\mathcal{O}_B$ is simple, 
whenever matrix $B$ is irreducible (i.e. a certain power
of $B$ is a strictly positive integer matrix). 
It was established in [Cuntz  \&  Krieger 1980]  \cite{CuKr1},  that $K_0(\mathcal{O}_B)\cong {\Bbb Z}^n / (I-B^t){\Bbb Z}^n$
and $K_1(\mathcal{O}_B)=Ker~(I-B^t)$, where $B^t$ is a transpose of the matrix $B$.
It is not difficult to see, that whenever $det~(I-B^t)\ne 0$, the $K_0(\mathcal{O}_B)$
is a finite abelian group and $K_1(\mathcal{O}_B)=0$.  The both groups are invariants of the
stable isomorphism class of the Cuntz-Krieger algebra.

%**************************************************************************
\section{Proof of theorem \ref{thm1}}
%***************************************************************************
Let $p$ be such, that $\mathcal{E}_{CM}$ has a good reduction at $\mathfrak{P}$;  
the corresponding local zeta function   $\zeta_p(\mathcal{E}_{CM},z)=(1-tr~(\psi_{E(K)}(\mathfrak{P}))z+pz^2)^{-1}$,
where $\psi_{E(K)}$ is the Gr\"ossencharacter on $K$ and  $tr$ is the trace of algebraic 
number.  We have to prove,  that $\zeta_p(\mathcal{E}_{CM},z)=\zeta_p(\mathcal{ A}_{RM},z):=(1-tr~(A^{\pi(p)})z+pz^2)^{-1}$;
the last equality is a consequence of definition of $\zeta_p(\mathcal{A}_{RM},z)$.
Let $\mathcal{E}_{CM}\cong {\Bbb C}/L_{CM}$,  where $L_{CM}={\Bbb Z}+{\Bbb Z}\tau$
is a lattice in the complex plane  [Silverman 1994]  \cite[pp. 95-96]{S};  
let $K_0(\mathcal{A}_{RM})\cong\Lambda_{RM}$,  where 
$\Lambda_{RM}={\Bbb Z}+{\Bbb Z}\theta$ is a  pseudo-lattice in 
${\Bbb R}$. 
Roughly speaking,  we construct an invertible element (a unit) $u$ of
the ring $End~(\Lambda_{RM})$ attached to pseudo-lattice $\Lambda_{RM}=F(L_{CM})$,
  such  that:  
%***********************************************************
\begin{equation}
tr~(\psi_{E(K)}(\mathfrak{P}))=tr~(u)=tr~(A^{\pi(p)}).
\end{equation}
%***********************************************************  
The latter  will be achieved with the help of an explicit formula 
connecting endomorphisms of lattice $L_{CM}$ with such of 
the pseudo-lattice $\Lambda_{RM}$ \cite[p. 142]{Nik2}:
%***********************************************************
\begin{equation}\label{eq2}
\left(\begin{matrix}a & b\cr c & d\end{matrix}\right)\in End~(L_{CM})
\longmapsto
 \left(\begin{matrix}a & b\cr -c & -d\end{matrix}\right)\in End~(\Lambda_{RM}).
\end{equation}
%***********************************************************

\bigskip
We shall split the proof into a series of lemmas,
 starting with the following simple
%*********************************************************
\begin{lemma}\label{lm1}
Let $A=(a,b,c,d)$ be an integer matrix with $ad-bc\ne 0$
and $b=1$. Then $A$ is similar to the matrix
$(a+d, 1, c-ad, 0)$.
\end{lemma}
%********************************************************
\begin{proof}
Indeed, take a matrix $(1,0,d,1)\in SL_2({\Bbb Z})$.
The matrix realizes the similarity, i.e.
%***********************************************************
\begin{equation}\label{eq3}
\left(\begin{matrix} 1 & 0\cr -d & 1\end{matrix}\right)
\left(\begin{matrix} a & 1\cr  c & d\end{matrix}\right)
\left(\begin{matrix} 1 & 0\cr d & 1\end{matrix}\right)=
\left(\begin{matrix} a+d & 1\cr c-ad & 0\end{matrix}\right).
\end{equation}
%***********************************************************  
\end{proof}
%*********************************************************
\begin{lemma}\label{lm2}
The matrix $A=(a+d, 1, c-ad, 0)$ is similar to its
transpose $A^t=(a+d, c-ad, 1, 0)$. 
\end{lemma}
%********************************************************
\begin{proof}
We shall use the following criterion: the
(integer) matrices $A$ and $B$ are similar, if and only if 
the characteristic matrices $xI-A$ and $xI-B$ have the same Smith normal
form.  The calculation for the matrix $xI-A$ gives:
%***************************************************************
$$
\left(\begin{matrix} x-a-d & -1\cr ad-c & x\end{matrix}\right)\sim
\left(\begin{matrix} x-a-d & -1\cr  x^2-(a+d)x+ad-c & 0\end{matrix}\right)\sim
$$
$$
\sim \left(\begin{matrix} 1 & 0\cr 0 & x^2-(a+d)x+ad-c\end{matrix}\right),
$$
%***********************************************************  
where $\sim$ are the elementary operations between the rows (columns)  
of the matrix.  Similarly, a calculation for the matrix $xI-A^t$
gives:
%***************************************************************
$$
\left(\begin{matrix} x-a-d & ad-c\cr -1 & x\end{matrix}\right)\sim
\left(\begin{matrix} x-a-d & x^2-(a+d)x+ad-c\cr -1& 0\end{matrix}\right)\sim
$$
$$
\sim\left(\begin{matrix} 1 & 0\cr 0 & x^2-(a+d)x+ad-c\end{matrix}\right).
$$
%***********************************************************
Thus, $(xI-A)\sim (xI-A^t)$ and lemma \ref{lm2} follows.
\end{proof}
%*********************************************************
\begin{corollary}\label{cr1}
The matrices $(a, 1, c, d)$ and $(a+d, c-ad, 1, 0)$
are similar. 
\end{corollary}
%********************************************************
\begin{proof}
 It follows from lemmas \ref{lm1}-\ref{lm2}.
\end{proof}

\medskip
Let $\mathcal{E}_{CM}$ be elliptic curve with the complex multiplication
by an order $R$ in the ring of integers of the imaginary quadratic
field $k$. Then $\mathcal{A}_{RM}=F(\mathcal{E}_{CM})$ is a noncommutative torus
with real multiplication by the order $\mathfrak{R}_n$ of conductor $n\ge 1$  in the ring of
integers $O_\mathfrak{k}$  of a  real quadratic field $\mathfrak{k}$. Let $tr~(\alpha)=\alpha+\bar\alpha$ be the
trace function of a (quadratic) algebraic number field. 
%*********************************************************
\begin{lemma}\label{lm3}
Each $\alpha\in R$ goes under $F$ into an $\omega\in \mathfrak{R}$,
such that $tr~(\alpha)=tr~(\omega)$. 
\end{lemma}
%********************************************************
\begin{proof}
Recall that each $\alpha\in R$ can be written  in a matrix
form for a given base $\{\omega_1,\omega_2\}$ of the lattice
$\Lambda$. Namely,
 %*******************************************************************
\begin{equation}\label{eq4}
\left\{
\begin{array}{cc}
\alpha\omega_1 &= a\omega_1 +b\omega_2\\
\alpha\omega_2 &= c\omega_1 +d\omega_2,
\end{array}
\right.
\end{equation}
%********************************************************* 
where $(a,b,c,d)$ is an integer matrix with $ad-bc\ne 0$. 
and $tr~(\alpha)=a+d$.

The first equation implies $\alpha=a+b\tau$;  since both $\alpha$ and
$\tau$ are algebraic integers, one concludes that $b=1$.

In view of corollary \ref{cr1}, in a base $\{\omega_1',\omega_2'\}$,
the $\alpha$ has a matrix form $(a+d, c-ad, 1, 0)$. 

To calculate a real quadratic $\omega\in \mathfrak{R}$ corresponding to $\alpha$,
recall an explicit formula from \cite[p.142]{Nik2}. Namely,  
every endomorphism $(a,b,c,d)$ of the lattice $L_{CM}$ maps to the endomorphism
$(a,b,-c,-d)$ of the pseudo-lattice $\Lambda_{RM}=F(L_{CM})$.
Thus, one gets a map:
%***********************************************************
\begin{equation}\label{eq5}
\left(\begin{matrix} a+d & c-ad\cr 1 & 0\end{matrix}\right)
\longmapsto
 \left(\begin{matrix} a+d & c-ad\cr -1 & 0\end{matrix}\right). 
\end{equation}
%***********************************************************
In other words, for a given base $\{\lambda_1,\lambda_2\}$ of
the pseudo-lattice ${\Bbb Z}+{\Bbb Z}\theta$ one can write
%*******************************************************************
\begin{equation}\label{eq6}
\left\{
\begin{array}{cc}
\omega\lambda_1 &= (a+d)\lambda_1 +(c-ad)\lambda_2\\
\omega\lambda_2 &= -\lambda_1.
\end{array}
\right.
\end{equation}
%********************************************************* 
It is an easy exercise to verify that $\omega$ is a real 
quadratic integer with $tr~(\omega)=a+d$. The latter coincides
with the $tr~(\alpha)$.
\end{proof}

\medskip
Let $\omega\in \mathfrak{R}$ be an endomorphism of the pseudo-lattice
$\Lambda={\Bbb Z}+{\Bbb Z}\theta$ of degree $deg~(\omega):=\omega\bar\omega=n$.
The endomorphism maps $\Lambda$  to a sub-lattice, $\Lambda_n$,  of index $n$;
any such has the form $\Lambda_n={\Bbb Z}+(n\theta){\Bbb Z}$ [Borevich \& Shafarevich 1966]  \cite[p.131]{BS}. 
Moreover,  $\omega$ generates an automorphism, $u$,  of the pseudo-lattice
$\Lambda_n$;  the two are related by the following lemma.  
%*********************************************************
\begin{lemma}\label{lm4}
$tr~(u)=tr~(\omega)$. 
\end{lemma}
%********************************************************
\begin{proof}
Let us calculate the action of endomorphism 
$\omega=(a+d, c-ad, -1, 0)$ on the pseudo-lattice 
$\Lambda_n={\Bbb Z}+(n\theta){\Bbb Z}$.
Since $deg~(\omega)=c-ad=n$,  one gets
%*******************************************************************
\begin{equation}\label{eq7}
\left(\begin{matrix}a+d & n\cr -1 & 0\end{matrix}\right)
\left(\begin{matrix} 1\cr \theta\end{matrix}\right)=
\left(\begin{matrix} a+d & 1\cr -1 & 0\end{matrix}\right)
\left(\begin{matrix} 1\cr n\theta\end{matrix}\right),
\end{equation}
%********************************************************* 
where $\{1,\theta\}$ and $\{1,n\theta\}$ are bases of the pseudo-lattices 
$\Lambda$ and $\Lambda_n$, respectively,  and $u=(a+d, 1, -1, 0)$
is an automorphism of $\Lambda_n$.  It is easy to see,  that
$tr~(u)=a+d=tr~(\omega)$.   
Lemma \ref{lm4} follows.
\end{proof}

%**************************************************************************
\begin{remark}\label{rm1}
There exists a canonical proof of lemma \ref{lm4} based on the notion of
a subshift of finite type [Wagoner 1999]  \cite{Wag1};  we shall give such a proof  below,  since
it generalizes to   pseudo-lattices of any rank.   
Consider a dimension group ([Blackadar 1986]  \cite[p.55]{B}) corresponding to the endomorphism $\omega$
of lattice ${\Bbb Z}^2$,  i.e. the limit $G(\omega)$: 
%************************************************************************************
\begin{equation}
{\Bbb Z}^2 \buildrel\omega \over\to 
{\Bbb Z}^2 \buildrel\omega \over\to 
{\Bbb Z}^2 \buildrel\omega \over\to
\dots
\end{equation}
%*********************************************************************************** 
 It is known that $G(\omega)\cong {\Bbb Z}[{1\over\lambda}]$, where $\lambda>1$ is the Perron-Frobenius
 eigenvalue of $\omega$.  We shall write $\hat\omega$ to denote the shift automorphism of dimension
 group $G(\omega)$, ([Effros  1981]  \cite[p.37]{E})  and 
 $\zeta_{\omega}(t)=\exp\left(\sum_{k=1}^{\infty}{tr~(\omega^k)\over k}t^k\right)$ and 
 $\zeta_{\hat\omega}(t)=\exp\left(\sum_{k=1}^{\infty}{tr~(\hat\omega^k)\over k}t^k\right)$ the 
 corresponding Artin-Mazur zeta functions \cite[p.273]{Wag1}.     
 Since the Artin-Mazur zeta function of the  subshift of finite type  is an invariant of shift equivalence,
 we conclude that  $\zeta_{\omega}(t)\equiv\zeta_{\hat\omega}(t)$;  in particular, $tr~(\omega)=tr~(\hat\omega)$.  
 Hence the matrix form of $\hat\omega=(a+d, 1, -1, 0)=u$ and, therefore, $tr~(u)=tr~(\omega)$. 
Lemma \ref{lm4} follows. 
\end{remark}
%***********************************************************************
%*********************************************************
\begin{lemma}\label{lm5}
The automorphism $u$ is a unit of the ring $\mathfrak{R}_n=End~(\Lambda_n)$;
it is the fundamental unit of $\mathfrak{R}_n$,  whenever $n=p$ (a prime number)
and  $tr~(u)=tr~(\psi_{E(K)}(\mathfrak{P}))$. 
\end{lemma}
%********************************************************
\begin{proof}
(i)  Since $deg~(u)=1$,  the element $u$ is invertible and, therefore,
a unit of the ring $\mathfrak{R}_n$;  in general, unit  $u$ is not the fundamental
unit of $\mathfrak{R}_n$,  since it is possible that  $u=\varepsilon^a$,  where
$\varepsilon$ is another unit of $\mathfrak{R}_n$ and $a\ge 1$.  

\medskip
(ii)  When $n=p$ is a prime number,  then we let   $\psi_{E(K)}(\mathfrak{P})$
be the corresponding Gr\"ossencharacter on $K$ attached to an elliptic
curve $\mathcal{E}_{CM}\cong E(K)$,  see Section 2.1 for the notation. 
The Gr\"ossencharacter can be identified with a complex number $\alpha\in k$
of the imaginary quadratic field $k$ associated to the complex multiplication.

Let  $tr~(u)=tr~(\psi_{E(K)}(\mathfrak{P}))$ and suppose to the contrary, that
$u$ is not the fundamental unit of $\mathfrak{R}_p$,  i.e.  $u=\varepsilon^a$
for a unit $\varepsilon\in \mathfrak{R}_p$ and an integer $a\ge 1$. 
Then there exists a Gr\"ossencharacter $\psi^{\prime}_{E(K)}(\mathfrak{P})$,
such that 
%*****************************************************************************
\begin{equation}
tr~(\psi^{\prime}_{E(K)}(\mathfrak{P}))< tr~(\psi_{E(K)}(\mathfrak{P})). 
\end{equation}
%***************************************************************************
Since $tr~(\psi_{E(K)}(\mathfrak{P}))=q_\mathfrak{P}+1-\# \tilde E({\Bbb F}_\mathfrak{P})$,
one concludes that  $\# \tilde E({\Bbb F}^{\prime}_\mathfrak{P})  > \# \tilde E({\Bbb F}_\mathfrak{P})$;
in other words,  there exists a non-trivial extension 
${\Bbb F}^{\prime}_\mathfrak{P}\supset {\Bbb F}_\mathfrak{P}$
of the finite field  ${\Bbb F}_\mathfrak{P}$.  The latter is impossible,
since any extension of  ${\Bbb F}_\mathfrak{P}$ has the form ${\Bbb F}_{\mathfrak{P}^n}$
for some $n\ge 1$;  thus $a=1$, i.e. unit $u$ is the fundamental unit 
of the ring $\mathfrak{R}_p$. Lemma \ref{lm5} follows. 
\end{proof}

%*********************************************************
\begin{lemma}\label{lm6}
$tr~(\psi_{E(K)}(\mathfrak{P}))=tr~(A^{\pi(p)})$. 
\end{lemma}
%********************************************************
\begin{proof}
In view of lemma \ref{lm10},  the fundamental
unit of the order $\mathfrak{R}_p$ is given by the formula $\varepsilon_p=\varepsilon^{\pi(p)}$,
where $\varepsilon$ is the fundamental unit of the ring $O_\mathfrak{k}$
and $\pi(p)$ an integer number. 
On the other hand,  matrix $A=\prod_{i=1}^n(a_i, 1, 1, 0)$,
where $\theta=(\overline{a_1,\dots,a_n})$ is a purely periodic
continued fraction.  Therefore
%************************************************************
\begin{equation}
A\left(\begin{matrix} 1\cr\theta\end{matrix}\right)=
\varepsilon\left(\begin{matrix} 1\cr\theta\end{matrix}\right), 
\end{equation}
%**************************************************************
where $\varepsilon>1$ is the fundamental unit of the real
quadratic field $\mathfrak{k}={\Bbb Q}(\theta)$.  In other words,
$A$ is the matrix form of  the fundamental unit $\varepsilon$. 
Therefore  the matrix form of the fundamental unit  $\varepsilon_p=\varepsilon^{\pi(p)}$ 
of $\mathfrak{R}_p$ is given by matrix   $A^{\pi(p)}$.
We apply lemma \ref{lm5} and get  
%*************************************************************
\begin{equation}\label{eq10}
tr~(\psi_{E(K)}(\mathfrak{P}))=tr~(\varepsilon_p)=tr~(A^{\pi(p)}).
\end{equation}
%*************************************************************
Lemma \ref{lm6} follows.
\end{proof}

\bigskip
One can finish the proof of theorem \ref{thm1} by comparing the local $L$-series 
of the Hasse-Weil $L$-function for the $\mathcal{E}_{CM}$ with  that of the local zeta  for
the $\mathcal{A}_{RM}$. 
The local $L$-series for $\mathcal{E}_{CM}$ are
$L_\mathfrak{P}(E(K),T)=1-a_\mathfrak{P}T+q_\mathfrak{P}T^2$ if the $\mathcal{E}_{CM}$ has a good reduction at $\mathfrak{P}$
and $L_\mathfrak{P}(E(K),T)=1-\alpha T$ otherwise; here 
%***************************************************************
\begin{eqnarray}\label{eq11}
 q_\mathfrak{P} &=& N^K_{\Bbb Q}\mathfrak{P}=\# {\Bbb F}_\mathfrak{P}=p,\cr
a_\mathfrak{P} &=& q_\mathfrak{P}+1-\# \tilde E({\Bbb F}_\mathfrak{P})=tr~(\psi_{E(K)}(\mathfrak{P})),\cr
\alpha &\in & \{-1,0,1\}. 
\end{eqnarray}
%*********************************************************************  
Therefore, 
%***********************************************************************
\begin{equation}
L_\mathfrak{P}(\mathcal{E}_{CM},T)=
\begin{cases}
1-tr~(\psi_{E(K)}(\mathfrak{P}))T+pT^2, & \mbox{for good reduction }  \cr
1-\alpha T,  & \mbox{for bad reduction}.
\end{cases}
\end{equation}
%*************************************************************************** 

\medskip
Let now $\mathcal{A}_{RM}=F(\mathcal{E}_{CM})$. 
%*******************************************************
\begin{lemma}\label{lm7}
$\zeta_p^{-1}(\mathcal{A}_{RM}, T)=1-tr~(A^{\pi(p)})T+pT^2$, whenever $p\nmid tr^2(A)-4$.
\end{lemma}
%*********************************************
{\it Proof.} By the formula   $K_0(\mathcal{O}_{B})={\Bbb Z}^2/(I-B^t){\Bbb Z}^2$,
one gets:
%*****************************************************************************
\begin{equation}
|K_0(\mathcal{O}_{L_p^n})|=\left|{{\Bbb Z}^2\over (I-(L_p^n)^t){\Bbb Z}^2}\right|=
|det(I-(L_p^n)^t)|=|Fix ~(L_p^n)|,
\end{equation}
%*****************************************************************
where $Fix~(L_p^n)$ is the set of (geometric) fixed points of the
endomorphism $L_p^n: {\Bbb Z}^2\to {\Bbb Z}^2$. Thus,
%***********************************************************************
\begin{equation}
\zeta_p(\mathcal{A}_{RM}, z)= \exp \left(\sum_{n=1}^{\infty}{|Fix~(L_p^n)|\over n} ~z^n\right), 
\quad z\in {\Bbb C}. 
\end{equation}
%*************************************************************************** 
But the latter series is an Artin-Mazur zeta function of the endomorphism $L_p$;
it converges to a rational function $det^{-1}(I-zL_p)$ [Hartshorn  1977]  \cite[p.455]{HAR}. Thus,
$\zeta_p(\mathcal{A}_{RM}, z)= det^{-1}(I-zL_p)$.

The substitution $L_p=(tr~(A^{\pi(p)}), p, -1, 0)$ gives us:
%******************************************************************
\begin{equation}
det~(I-zL_p)=
det~\left(\begin{matrix} 1-tr~(A^{\pi(p)})z & -pz\cr z & 1\end{matrix}\right)=1-tr~(A^{\pi(p)})z+pz^2.
\end{equation}
%***********************************************************
Put $z=T$ and get $\zeta_p(\mathcal{A}_{RM}, T)=(1-tr~(A^{\pi(p)})T+pT^2)^{-1}$,
which is a conclusion of lemma \ref{lm7}.
$\square$

%*******************************************************
\begin{lemma}\label{lm8}
$\zeta_p^{-1}(\mathcal{A}_{RM}, T)=1-\alpha T$, whenever $p ~| ~tr^2(A)-4$.
\end{lemma}
%*********************************************
\begin{proof}
Indeed, $K_0(\mathcal{O}_{1-\alpha^n})={\Bbb Z}/(1-1+\alpha^n){\Bbb Z}={\Bbb Z}/\alpha^n{\Bbb Z}$.
Thus, $|K_0(\mathcal{O}_{1-\alpha^n})|=det~(\alpha^n)=\alpha^n$. By the definition,
%*******************************************************************************
\begin{equation}
\zeta_p(\mathcal{A}_{RM},z)=\exp\left(\sum_{n=1}^{\infty} {\alpha^n\over n}z^n\right)=
\exp\left(\sum_{n=1}^{\infty} {(\alpha z)^n\over n}\right)=
{1\over 1-\alpha z}.
\end{equation}
%**********************************************************************
The substitution $z=T$ gives a conclusion of lemma \ref{lm8}.
\end{proof}

%*******************************************************
\begin{lemma}\label{lm9}
Let $\mathfrak{P}\subset K$ be a prime ideal over $p$; 
then $\mathcal{E}_{CM}=E(K)$ has a bad reduction at $\mathfrak{P}$
if and only if  $p~|~tr^2(A)-4$. 
\end{lemma}
%*********************************************
\begin{proof}
Let $k$ be a  field of complex 
multiplication of the $\mathcal{E}_{CM}$; its discriminant we
shall write as $\Delta_k<0$. It is known, that whenever $p~|~\Delta_k$,
the $\mathcal{E}_{CM}$ has a bad reduction at the prime ideal $\mathfrak{P}$ over $p$.

On the other hand, the explicit formula (\ref{eq2}) applied to the matrix $L_p$
gives us $F: (tr~(A^{\pi(p)}),  p,  -1, 0) \mapsto (tr~(A^{\pi(p)}), p,  1, 0)$.
The characteristic polynomials of the above matrices
are  $x^2-tr~(A^{\pi(p)})x+p$ and $x^2-tr~(A^{\pi(p)})x-p$, respectively.   
They generate an imaginary (resp., a real) quadratic field 
$k$ (resp., $\mathfrak{k}$) with the discriminant $\Delta_k=tr^2 (A^{\pi(p)})-4p<0$
(resp., $\Delta_\mathfrak{k}=tr^2 (A^{\pi(p)})+4p>0$). Thus,
$\Delta_\mathfrak{k}-\Delta_k=8p$. It is easy to see, that 
$p~|~\Delta_\mathfrak{k}$ if and only if  $p~|~\Delta_k$.
It remains to express the discriminant $\Delta_\mathfrak{k}$ in terms
of the matrix $A$. Since the characteristic polynomial for $A$
is $x^2-tr~(A)x+1$, it follows  that $\Delta_\mathfrak{k}=tr^2(A)-4$.
\end{proof}

\medskip
We are prepared now to prove the first part of theorem \ref{thm1}.
Note, that a critical piece of information is   provided by 
lemma \ref{lm6}, which says that $tr~(\psi_{E(K)}(\mathfrak{P}))=tr~(A^{\pi(p)})$.
Thus, in view of lemmas \ref{lm7}-\ref{lm9},  $L_\mathfrak{P}(\mathcal{E}_{CM},T)\equiv \zeta_p^{-1}(\mathcal{A}_{RM}, T)$.
The first part of theorem \ref{thm1} follows.

\subsection{Case  $p$ is a good prime}
 Let us prove the second part of theorem \ref{thm1} in the case
$n=1$. From the left side: $K_0(\mathcal{A}_{RM}\rtimes_{L_p}{\Bbb Z})\cong K_0(\mathcal{O}_{L_p})\cong
{\Bbb Z}^2/(I-L_p^t){\Bbb Z}^2$, where $L_p=(tr~(A^{\pi(p)}), p, -1, 0)$. To
calculate the abelian group ${\Bbb Z}^2/(I-L_p^t){\Bbb Z}^2$, we shall 
use a reduction of the matrix $I-L_p^t$ to the Smith normal form:
%******************************************************************
$$
I-L_p^t=
\left(\begin{matrix} 1-tr~(A^{\pi(p)}) & 1\cr -p & 1\end{matrix}\right)\sim
\left(\begin{matrix} 1+p-tr~(A^{\pi(p)}) & 0\cr -p & 1\end{matrix}\right)\sim
$$
$$
\sim\left(\begin{matrix} 1 & 0\cr 0 & 1+p-tr~(A^{\pi(p)})\end{matrix}\right).  
$$
%***************************************************************
Therefore, $K_0(\mathcal{O}_{L_p})\cong {\Bbb Z}_{1+p-tr~(A^{\pi(p)})}$. 

From the right side, the $\mathcal{E}_{CM}({\Bbb F}_\mathfrak{P})$ is an elliptic
curve over the field of characteristic $p$. Recall,  that the chord and tangent law
turns the $\mathcal{E}_{CM}({\Bbb F}_\mathfrak{P})$ into a finite abelian group. The group
is cyclic and has the order $1+q_\mathfrak{P}-a_\mathfrak{P}$.
But $q_\mathfrak{P}=p$ and $a_\mathfrak{P}= tr~(\psi_{E(K)}(\mathfrak{P}))=tr~(A^{\pi(p)})$
(lemma \ref{lm6}).  Thus, $\mathcal{E}_{CM}({\Bbb F}_\mathfrak{P})\cong {\Bbb Z}_{1+p-tr~(A^{\pi(p)})}$;
therefore  $K_0(\mathcal{O}_{L_p})\cong \mathcal{E}_{CM}({\Bbb F}_p)$.

The general case $n\ge 1$ is treated likewise. Repeating  the argument of lemmas \ref{lm1}-\ref{lm2},
it follows that $L_p^n=(tr~(A^{n\pi(p)}, p^n, -1, 0)$. Then one gets  $K_0(\mathcal{O}_{L_p^n})\cong {\Bbb Z}_{1+p^n-tr~(A^{n\pi(p)})}$
on the left side. From the right side, $|\mathcal{E}_{CM}({\Bbb F}_{p^n})|=1+p^n-tr~(\psi^n_{E(K)}(\mathfrak{P}))$;
but a repetition of the argument of lemma \ref{lm6} yields us $tr~(\psi^n_{E(K)}(\mathfrak{P}))=tr~(A^{n\pi(p)})$.
Comparing the left and right sides, one gets that $K_0(\mathcal{O}_{L_p^n})\cong \mathcal{E}_{CM}({\Bbb F}_{p^n})$.
This argument finishes the proof of the second part of theorem \ref{thm1} for the good primes.

\subsection{Case $p$ is a bad prime}
From the proof of lemma \ref{lm8}, one gets for the left side  $K_0(\mathcal{O}_{\varepsilon_n})\cong {\Bbb Z}_{\alpha^n}$.
From the right side, it holds    $|\mathcal{E}_{CM}({\Bbb F}_{p^n})|=1+q_\mathfrak{P}-a_\mathfrak{P}$,
where $q_\mathfrak{P}=0$ and $a_\mathfrak{P}=tr~(\varepsilon_n)=\varepsilon_n$. 
Thus,   $|\mathcal{E}_{CM}({\Bbb F}_{p^n})|=1-\varepsilon_n=1-(1-\alpha^n)=\alpha^n$. 
Comparing the left and right sides, we conclude that
$K_0(\mathcal{O}_{\varepsilon_n})\cong \mathcal{E}_{CM}({\Bbb F}_{p^n})$ at the bad primes.

\medskip
Since all  cases are exhausted,  theorem \ref{thm1} is proved.

\bibliographystyle{amsplain}

%**********************************************************

\end{document}